\documentclass[12pt]{amsart}

\usepackage{amssymb}
\usepackage{epsf}

\newtheorem{thm}{Theorem}[section]

\newtheorem{prop}[thm]{Proposition}
\newtheorem{conj}[thm]{Conjecture} 
   \newtheorem{example}[thm]{Example}
 
   \newtheorem{definition}[thm]{Definition}
\newenvironment{defn}{\begin{definition}\rm}{\end{definition}} 

\begin{document}

\title[Numerical Schubert Calculus]{Numerical Schubert Calculus}

\author{Birkett Huber, Frank Sottile \and Bernd Sturmfels}
\address{Department of Mathematics, Texas A\&M University}
\email[Birkett Huber]{birk@@math.cornell.edu}

\address{Department of Mathematics, University of Toronto,
100 St.~George Street, Toronto. Ont., M5S 3G3, CANADA}
\email[Frank Sottile]{sottile@@math.toronto.edu}

\address{Department of Mathematics, University of California,
 Berkeley CA 94720, USA}
\email[Bernd Sturmfels]{bernd@@math.berkeley.edu}
\thanks{first author supported in part by NSF grant DMS-9508742}
\thanks{second author supported in part by NSF grant DMS-9022140}
\thanks{third author supported in part by a David and Lucille Packard
	Fellowship and an NSF National Young Investigator Fellowship}
\date{9 June 1997}


\begin{abstract}
We develop numerical homotopy algorithms for solving systems 
of polynomial equations arising from the classical Schubert calculus.
These homotopies are  optimal in that generically no paths diverge.
For problems defined by hypersurface Schubert conditions we give 
two algorithms based on extrinsic deformations of the Grassmannian:
one is derived from a Gr\"obner basis for the Pl\"ucker ideal of the 
Grassmannian and the other from a SAGBI basis for its projective 
coordinate ring.  The more general case of special Schubert conditions is 
solved by delicate intrinsic deformations, called Pieri homotopies, which
first arose in the study of enumerative geometry over the real numbers.
Computational results are presented and
applications to control theory are discussed.
\end{abstract}

\maketitle
\tableofcontents

\section{Introduction}
Suppose we are given linear subspaces $K_1,\ldots,K_n$ of ${\bf C}^{m+p}$
with $\dim K_i = m+1-k_i$ and $k_1+\cdots+k_n = mp $.
Our problem is to find all $p$-dimensional linear subspaces of 
${\bf C}^{m+p}$ which meet each $K_i$ nontrivially.
When the given linear subspaces are in general position, the condition
$k_1+\cdots+k_n=mp$ guarantees that there is a finite number
$\,d \,=\, d(m,p,k_1,\ldots,k_n)\,$ of such $p$-planes.
The classical Schubert calculus~\cite{Kleiman_Laksov} gives 
the following recipe for computing the number $d$.
Let $h_1,\ldots,h_m$ be indeterminates with $degree(h_i) = i$.
For each integer sequence $\lambda_1\geq \cdots\geq \lambda_{p+1} $
we define the following polynomial:
\begin{equation}\label{schur}
S_\lambda\ :=\ \det ( h_{\lambda_i+j-i})_{1\leq i,j\leq p+1}.
\end{equation}
Here $h_0 := 1$ and $h_i:=0$ if $i<0$ or $i>m$.
Let $I$ be the ideal in ${\bf Q}[h_1,\ldots,h_m]$ generated by those 
$S_\lambda$ with  $m \geq \lambda_1 $ and $\lambda_{p+1} \geq 1$.
The quotient ring ${\mathcal A}_{m,p}:= {\bf Q}[h_1,\ldots,h_m]/I$
is the cohomology ring of the Grassmannian
of $p$-planes in ${\bf C}^{m+p}$. It
is Artinian with one-dimensional socle in degree $mp$.
In the socle we have the relation
\begin{equation}
\label{soclerel}
 d \cdot (h_m)^p \, -  \,
h_{k_1} h_{k_2} \cdots h_{k_n}  \quad \in \quad I. 
\end{equation}
Thus we can compute the number $d$ by normal form reduction modulo 
any Gr\"obner basis for $I$. More efficient
methods for computing in the ring ${\mathcal A}_{m,p}$
are implemented in  the Maple package SF~\cite{Stembridge_SF}.

In the important special case $k_1=\cdots=k_n=1$
there is an explicit formula for $d$:
\begin{equation}\label{grassdeg} d \quad = \quad
 \frac{1! \, 2! \, 3! \cdots (p\!- \!2) ! \, (p \!-\!1)! \cdot
(mp)!}{m!\, (m \! + \! 1)! \, (m \! + \, 2)! 
\cdots(m \! + \! p \! - \! 1)!} . 
\end{equation}
The integer on the right hand side is the degree of the Grassmannian in
its Pl\"ucker embedding. This
formula is due to~\cite{Schubert_degree}; 
see also~\cite[XIV.7.8]{Hodge_Pedoe} and Section 2.3 below.

The objective of this paper is to present semi-numerical
algorithms for computing all $d$ solution planes from the input
data $K_1,\ldots,K_n$. This amounts to solving certain
systems of polynomial equations.  Our algorithms are based on 
the paradigm of {\it numerical homotopy
methods}~\cite{Morgan_book,Allgower_Georg,Allgower_Georg_handbook}. 

Homotopy methods have been developed for the following classes of
polynomial systems: 

\begin{enumerate}
\item complete intersections in affine or projective
spaces~\cite{Drexler,Garcia_Zangwill}, 
\item complete intersections in products of projective
spaces~\cite{MS_m-homogeneous}, 
\item complete intersections in toric
varieties~\cite{CVVerschelde,Huber_Sturmfels}.
\end{enumerate}

In these cases the number of paths to be traced is optimal and
equal to the standard combinatorial bounds:
\begin{enumerate}
\item the B\'ezout number (= the product of the degrees of the equations)
\item the generalized B\'ezout number  for multihomogeneous systems
\item the BKK bound~\cite{Bernstein,Kouchnirenko,Khovanskii_newton}
(= mixed volume of the Newton polytopes)
\end{enumerate}

None of these known homotopy methods is applicable to our problem, as the
following simple example shows: Take $m=3$, $p=2$, and $k_1=\cdots=k_6=1$,
that is, we seek the $2$-planes in ${\bf C}^5$ which meet six
general $3$-planes nontrivially. By formula (1.3) there are $d=5$ solutions.
Formulating this in Pl\"ucker coordinates gives $11$ 
homogeneous equations in ten variables,
the five quadrics in display~(\ref{I23})
below and six linear equations~(\ref{sixlin}).
A formulation in local coordinates~(\ref{sagbieqs})
has 6 quadratic equations in 6 unknowns, giving a B\'ezout bound of 64.
These 6 equations all have the same Newton polytope, which has 
normalized volume 17, giving a BKK bound of 17.

In Section 2 we give two homotopy algorithms
which solve our problem in the special case
$k_1 = \cdots =k_n = 1$, when the number of solutions equals (1.3).
The first algorithm is derived from a Gr\"obner basis for the Pl\"ucker 
ideal of a Grassmannian and the second from a SAGBI basis
for its projective coordinate ring.  (See~\cite{CHV}
or~\cite[Ch.~11]{Sturmfels_GPCP} for an introduction to SAGBI bases).
Both the {\it Gr\"obner homotopy} and {\it SAGBI homotopy} 
are techniques for finding linear sections of Grassmannians in their 
Pl\"ucker embedding. 

In Section 3 we address the general case of our problem, that is, 
we describe a numerical method for solving the polynomial equations
defined by {\it special Schubert conditions}. This is accomplished 
by applying a sequence of delicate intrinsic deformations, 
called {\it Pieri homotopies}, which were introduced
in~\cite{Sottile_explicit_pieri}. 
Pieri homotopies first arose in the study of enumerative geometry over the 
real numbers~\cite{Sottile_real_lines,Sottile_santa_cruz}. 
For the experts we remark that it is an open problem  to find 
{\it Littlewood-Richardson homotopies},
which would be relevant for solving polynomial equations defined by 
general Schubert conditions. 

A main challenge in designing homotopies for the Schubert calculus
is that one is not dealing with complete intersections: there are
generally more equations than variables. In Section 4 we discuss
some of the numerical issues arising from this challenge,
and how we propose to resolve them.
In Section 5 we discuss applications of these algorithms to control
theory and real enumerative geometry.
Finally, in Section 6 we present computational results. 

In closing the introduction
 let us emphasize that all homotopies described in this paper
are optimal in the sense the  number of  paths to be traced equals the 
number $d$. 
This means that for generic input data $K_1,\ldots,K_n$ no paths diverge.

\section{Linear equations in Pl\"ucker coordinates}
The set of $p$-planes in ${\bf C}^{m+p}$, $\mbox{\em Grass}(p,m+p)$, is
called the {\em Grassmannian of $p$-planes in  ${\bf C}^{m+p}$}.
This complex manifold of dimension $mp$ is naturally a subvariety
of the complex projective space ${\bf P}^{\binom{m+p}{p}-1} $.
To see this, represent a $p$-plane in ${\bf C}^{m+p}$ as the
column space of an $(m+p)\times p$-matrix $ X = (x_{ij})$.
The {\it Pl\"ucker coordinates} of that $p$-plane are the 
maximal minors of $X$, indexed by the set $\binom{[m+p]}{p}$
of sequences $\alpha: 1\leq \alpha_1<\alpha_2<\cdots<\alpha_p\leq m+p$:
\begin{equation}\label{minors}
 [\alpha]\ \longrightarrow\  \det
\left[\begin{array}{ccc}
x_{\alpha_1\,1}&\cdots& x_{\alpha_1\,p}\\
\vdots&\ddots&\vdots\\
x_{\alpha_p\,1}&\cdots& x_{\alpha_p\,p}\end{array}\right] .
\end{equation}

This section deals with the ``$k_i = 1$'' case  of
the problem stated in the Introduction.
Given $mp$ general $m$-planes $K_1,\ldots,K_{mp}$, we 
wish to find all $p$-planes $X$ which meet $K_1,\ldots,K_{mp}$ nontrivially.
This geometric condition translates into linear equations in the 
Pl\"ucker coordinates:
Represent $X$ as an $(m+p) \times p$-matrix as above,
represent $K_i$ as an $(m+p)\times m$-matrix, and form the
$(m+p) \times (m+p)$-matrix $ [\, X \mid K_i \,]  $. Then
$$ \,X \,\cap \,K_i \,\not= \,\{ 0 \} \qquad
\hbox{if and only if} \qquad \det \,[\,X \mid K_i \,] \,\,\,= 
\,\,\, 0 . $$
Laplace expansion with respect to the first $p$ columns gives
\begin{equation}\label{lineqs}
 \det \, [\, X \mid K_i \,] \quad = \quad
  \sum_{\alpha\in\binom{[m+p]}{p}} 
C_{\alpha}^{i} \cdot [\alpha], 
\end{equation}
where $C_{\alpha}^{i}$ is the correctly signed
maximal minor of $K_i$ complementary to $\alpha$.
Hence our problem is to solve $mp$ linear
equations (\ref{lineqs}) on the Grassmannian.
The number of solutions is the degree of the Grassmannian in its
Pl\"ucker embedding, which is given in (1.3).

The Grassmannian is represented either {\it implicitly}, 
as the zero set of  polynomials in the  Pl\"ucker coordinates, or 
{\it parametrically}, as the image of the polynomial map (\ref{minors}).
These two representations lead to two different numerical homotopies.
The implicit representation gives
the Gr\"obner homotopy in Section~\ref{grobnerhomotopy}
while the parametric representation gives the SAGBI homotopy in
Section~\ref{SAGBIhomotopy}. 
The first is conceptually simpler but the second
 is more efficient. In both methods the
number of paths to be traced equals the optimal number in (\ref{grassdeg}).

\subsection{An example}

We describe the two approaches for the case $(m,p)=(3,2)$.
The Grassmannian of 2-planes in ${\bf C}^5$ has dimension 6 and is 
embedded into ${\bf P}^9$. Its degree (\ref{grassdeg}) is five.
 The Gr\"obner homotopy works directly
in the ten Pl\"ucker coordinates:
$$ [12],\ [13],\ [14],\ [15],\ [23],\ [24],\ [25],\ [34],\ [35],\ [45]. $$
The ideal $I_{3,2}$ of the Grassmannian in the
Pl\"ucker embedding is generated by five quadrics:
\begin{equation}\label{I23}
\begin{array}{c}
\underline{[14][23]} \ -\ [13][24]\ +\ [12][34], \\ 
\underline{[15][23]}\ -\ [13][25]\ +\ [12][35], \\ 
\underline{[15][24]}\ -\ [14][25]\ +\ [12][45] , \\ 
\underline{[15][34]}\ -\ [14][35]\ +\ [13][45] , \\ 
\underline{[25][34]}\ -\ [24][35]\ +\ [23][45].
\end{array}
\end{equation}
This set is the reduced Gr\"obner basis for $I_{3,2}$ with
respect to any term order which selects the underlined terms
as leading terms (see Proposition 2.1 below).

Our problem is to compute all $2$-planes which meet
six sufficiently general $3$-planes $K_1,\ldots, K_6$ nontrivially.
This amounts to solving  (\ref{I23}) together with
six linear equations
\begin{equation}\label{sixlin}
\begin{array}{c}
 C_{12}^{i} \cdot [12] \, + \, 
C_{13}^{i} \cdot [13] \, + \, 
C_{14}^{i} \cdot [14] \, + \, 
C_{15}^{i} \cdot [15] \, + \, 
C_{23}^{i} \cdot [23]
\phantom{\, = \, 0}
\\
 \, + \,  C_{24}^{i} \cdot [24] \, + \, 
C_{25}^{i} \cdot [25] \, + \, 
C_{34}^{i} \cdot [34] \, + \, 
C_{35}^{i} \cdot [35] \, + \, 
C_{45}^{i} \cdot [45] 
\,\,\, = \,\,\, 0 ,
\end{array}
\end{equation}
for $i=1,\ldots,6 $.
This is an overdetermined system of
$11$ equations in $10$ homogeneous variables.
To solve it we introduce a parameter $t$ into (\ref{I23}) as follows:
$$ \begin{array}{l}
 [14][23]\ -\ t\;\,\cdot [13][24]\ +\ t^2 \cdot [12][34] \quad =\quad 0, \\ 
{[15][23]\ -\  t^2 \cdot [13][25]\ +\ t^4 \cdot [12][35]}\quad = \quad 0, \\ 
{[15][24]\ -\ t\;\,\cdot [14][25]\ +\ t^5 \cdot [12][45]}\quad =\quad 0, \\ 
{[15][34]\ -\  t^2 \cdot [14][35]\ +\ t^4 \cdot [13][45]}\quad = \quad 0, \\ 
{[25][34]\ -\ t\;\,\cdot [24][35]\ +\ t^2 \cdot [23][45] \quad =\quad 0.}
\end{array} \eqno (\ref{I23}') $$
We call (\ref{I23}$'$) the {\it Gr\"obner homotopy} because this is an
instance of the flat deformation which exists for any Gr\"obner basis;
see~\cite[Theorem 15.17]{Eisenbud_geometry}.  The flatness of this 
deformation ensures that, for almost every complex number $t$, the 
combined system $(\ref{I23}') \& (\ref{sixlin})$ has five roots. 

For $t\!=\!0$ the equations (\ref{I23}$'$) are square-free monomials.
We decompose their ideal:
\begin{equation}\label{primedec}
  \begin{array}{c}
 \langle \,\, [14][23],\,  [15][23],\,
 [15][24],\, [15][34],\,    [25][34] \,\, \rangle 
\hspace{2.2in}  \\\hspace{1.3in} = \quad  \,
\langle \, [23], [24], [34] \,\rangle  \,\, \cap \,\, 
\langle \, [15], [23], [34] \,\rangle  \,\, \cap \,\,
\langle \, [15], [23], [25] \,\rangle \\ 
\hspace{1.8in} \,\, \cap \,\, \,\,
\langle \, [14], [15], [34] \,\rangle  \,\, \, \cap \,\,\,
\langle \, [14], [15], [25] \,\rangle  .
\end{array}
\end{equation}
The five distinct solutions for $t=0$ are computed
by setting each listed triple of variables to zero and
then solving the six linear equations (\ref{sixlin}) in the remaining
seven variables. Thereafter we trace the
five solutions from $t=0$ to $t=1$ by numerical path continuation.
At $t=1$ we get the five solutions to our original problem.

We next describe the SAGBI homotopy. For this we choose
the local coordinates 
$$ X \quad = \quad
\left[\begin{array}{ccc}
            1   &    0   \\
         x_{21} & x_{22} \\
         x_{31} & x_{32} \\
         x_{41} & x_{42} \\
            0   &    1 
\end{array}
\right] 
$$
on the Grassmannian.
Substituting $X$ into (\ref{sixlin}) we get six polynomials in six unknowns:
\begin{equation}\label{sagbieqs}
 \begin{array}{ccc}
 C_{23}^{i} \cdot (x_{21} x_{32} -  x_{22} x_{31}) 
\,+\,   C_{24}^{i} \cdot (x_{21} x_{42} -  x_{22} x_{41})
\,+\,   C_{34}^{i} \cdot (x_{31} x_{42} -  x_{32} x_{41}) \\
\,+\,   C_{25}^{i} \cdot x_{21}
+   C_{35}^{i} \cdot x_{31}
+   C_{45}^{i} \cdot x_{41}
+    C_{12}^{i} \cdot  x_{22}
+   C_{13}^{i} \cdot  x_{32}
+  C_{14}^{i} \cdot  x_{42}
+   C_{15}^{i}.  \\
\end{array} 
\end{equation}
To solve these six equations
we introduce a parameter $t$ as follows:
$$ \begin{array}{ccc}
 C_{23}^{i} (x_{21} x_{32} - t \cdot x_{22} x_{31}) 
+   C_{24}^{i} (x_{21} x_{42} - t^2 \cdot x_{22} x_{41})
+   C_{34}^{i} (x_{31} x_{42} - t \cdot x_{32} x_{41}) \!\!\! \\
+   C_{25}^{i} \cdot x_{21}
+   C_{35}^{i} \cdot x_{31}
+   C_{45}^{i} \cdot x_{41}
+    C_{12}^{i} \cdot  x_{22}
+   C_{13}^{i} \cdot  x_{32}
+  C_{14}^{i} \cdot  x_{42}
+   C_{15}^{i}.  \\
\end{array} \eqno (\ref{sagbieqs}')
$$
The system (\ref{sagbieqs}$'$) has five complex roots
for almost all $t \in {\bf C}$. For $t=0$ we get a 
{\it generic unmixed sparse system}
(in the sense of~\cite{Huber_Sturmfels}) with support
$$ {\mathcal A} \quad = \quad 
\{\, 1, \,x_{21} , \, x_{22} \,,\, x_{31} \,, \, x_{32} \,, 
x_{41}\,,   x_{42}\,,\,
x_{21} x_{32} \,, \, x_{21} x_{42}\,, \,x_{31} x_{42} \,\}.$$
We identify ${\mathcal A}$ with a set of ten points in ${\bf Z}^6$.
Their convex hull {\it conv}$({\mathcal A})$ is a $6$-dimensional 
polytope with normalized volume five.  
We can therefore solve $(\ref{sagbieqs}')$ for $t=0$ using the homotopy
method in~\cite{Huber_Sturmfels} or~\cite{CVVerschelde},
provided the input data $K_1,\ldots,K_6$ are sufficiently generic.
Tracing the five roots from $t=0$ to $t=1$ by numerical path continuation,
we obtain the five solutions to our original problem.

\subsection{Gr\"obner homotopy}\label{grobnerhomotopy}
\begin{figure}[htb]
$$
\epsfxsize=1.05in\epsfbox{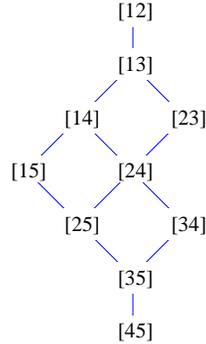}
$$
\caption{Young's poset for $(m,p) = (3,2)$.\label{youngs_poset}}
\end{figure}

We next describe a quadratic Gr\"obner basis
for the defining ideal of $\mbox{\em Grass}(p,m+p)$.
Let $S$ be the polynomial ring
over ${\bf C}$ in the variables $\,[\alpha]\,$
where $\alpha \in  \binom{ [m+p] }{p}$.
We define a partial order on these variables
as follows: $\,[\alpha] \leq [\beta]\,$ if and only if
$\alpha_i \leq \beta_i$ for $i=1,\ldots,p$.
This partially ordered set is called {\it Young's poset}.
Figure~\ref{youngs_poset} shows Young's poset for $(m,p) = (3,2)$.

Fix any linear ordering on the variables in $S$ which refines
the ordering in Young's poset, and let
 $\prec$ denote the induced degree reverse lexicographic term order on $S$.
Let $I_{m,p}$ be the ideal of polynomials  in $S$  which vanish
on the Grassmannian, that is, $I_{m,p}$ is the ideal of algebraic relations
among the maximal minors of a generic $(m+p) \times p$-matrix $X$.
The Gr\"obner homotopy is based on the following well-known result.

\begin{prop}\label{prop:inImp}
The initial ideal $\,in_\prec(I_{m,p}) \,$ is generated
by all quadratic monomials $\,[\,\alpha \,][\,\beta\,] \,$
where $\alpha_i < \beta_i$ and
$\alpha_j > \beta_j \,$  for some $i,j \in \{1,\ldots,p\}$.
\end{prop}

In other words, $\,in_\prec(I_{m,p}) \,$ is generated
by products of incomparable pairs in Young's poset.
Let ${\mathcal C}_{m,p}$ denote the set of all maximal chains
in  Young's poset.
For example,
\begin{eqnarray*}
{\mathcal C}_{3,2} &=&
\left\{ \, \{[12],[13],[14],[15],[25],[35],[45]\} \,,
\,\, \{[12],[13],[14],[24],[25],[35],[45]\}\,,\right. \\
&&\hspace{6.6pt}\{[12],[13],[14],[24],[34],[35],[45]\}\,,\,\,
 \{[12],[13],[23],[24],[34],[35],[45]\}\,,\\
&&\hspace{5.2pt}\left.  \{[12],[13],[23],[24],[25],[35],[45]\}\,\right\}
\end{eqnarray*}
A standard result in combinatorics~\cite{Stanley_schubert}  states
that the cardinality of ${\mathcal C}_{m,p}$ equals the number (\ref{grassdeg}). 
{}From Proposition~\ref{prop:inImp} we read off the following prime
decomposition which generalizes~(\ref{primedec}):
\begin{equation}\label{genprimedec}
 in_\prec(I_{m,p}) 
\quad = \quad
\bigcap_{C \in {\mathcal C}_{m,p}} \langle \,\,
[\alpha]\,\,: \,\, [\alpha] \not\in C \,\rangle . 
\end{equation}

For a proof of Proposition~\ref{prop:inImp} see~\cite[\S XIV.9]{Hodge_Pedoe}
or~\cite[Theorem (4.3)]{Bruns_Vetter} or~\cite[\S 3.1]{Sturmfels_invariant}.
In these references one finds
an explicit minimal Gr\"obner basis for $I_{m,p}$, which is
classically called the set of {\it straightening syzygies}. 
In the special case $p=2$ the straightening syzygies coincide with
the reduced Gr\"obner basis:

\begin{prop}
If $p=2$ then the reduced Gr\"obner basis of $I_{m,p}$
consists of the three-term  Pl\"ucker relations
$\,\underline{[i l][k j]} - [i k][j l] + [i j][k  l] \,$
where $1 \leq i < j < k < l \leq m+p = m+2$.
\end{prop}

For $p \geq 3$ the straightening syzygies and the
reduced Gr\"obner basis do not coincide, and they are complicated to describe.
For our purposes the following coarse description suffices.
Let  $Std$ be the set of all quadratic monomials in $S$
which do not lie in $in_\prec(I_{m,p})$.
The reduced Gr\"obner basis consists of elements of the form
\begin{equation}\label{redGB}
[\alpha][\beta] \quad - \quad
\sum_{[\gamma][\delta]  \in Std} 
E^{\alpha,\beta}_{\gamma,\delta} \cdot [\gamma][\delta] 
\end{equation}
where $[\alpha][\beta]$ runs over all generators
of $in_\prec(I_{m,p})$. The constants
$\,E^{\alpha,\beta}_{\gamma,\delta} \,$ are integers which
can be computed  by substituting (\ref{minors}) into (\ref{redGB}) and
solving linear equations.

The term order $\prec$ can be realized for the ideal $I_{m,p}$
by the following choices of weights. We define the {\it weight} of
the variable $[\alpha] = [\alpha_1 \alpha_2 \cdots \alpha_p]$
to be
\begin{equation}\label{weights}
 v_\alpha \quad := \quad - \frac{1}{2} \sum_{1 \leq i < j \leq p}
( \alpha_j - \alpha_i - 1)^2 . 
\end{equation}
If we replace each variable $\,[\alpha]\,$ in (\ref{redGB}) by
$\,[\alpha] \cdot t^{v_\alpha}\,$ and clear $t$-denominators afterwards,
then we get  the {\it Gr\"obner homotopy}:
$$ 
[\alpha][\beta] \quad - \quad
\sum_{[\gamma][\delta]  \in Std} 
E^{\alpha,\beta}_{\gamma,\delta} \cdot [\gamma][\delta] \cdot
t^{v_\gamma + v_\delta - v_\alpha - v_\beta}
\eqno (\ref{redGB}') 
$$
It can be checked that all exponents
$v_\gamma + v_\delta - v_\alpha - v_\beta$ appearing here
are positive integers. The special case $(m,p)=(3,2)$
is presented in (\ref{I23}$'$).

In the Gr{\"o}bner homotopy algorithm, we first solve
systems of $mp$ linear equations, one for each chain 
$C \in {\mathcal C}_{m,p}$.
These systems consist of the $mp$ equations~(\ref{lineqs}), one for each
$K_i$, and the  $\binom{m+p}{p}-mp-1$ equations
$$
[\alpha]\qquad\mbox{for}\qquad [\alpha]\not\in C,
$$ 
suggested by the prime decomposition~(\ref{genprimedec}).
Once this is accomplished, we 
trace each of these $d$ solutions from $t=0$ to $t=1$ in the
Gr\"obner homotopy~(\ref{I23}$'$).

Clearly, the weights $v_\alpha$ of~(\ref{weights}) are not best
possible for any  specific value of $m$ and $p$.
Smaller weights can be found using Linear Programming,
as explained e.g.~in the proof of~\cite[Proposition 1.11]{Sturmfels_GPCP}.
Another method would be to adapt the ``dynamic'' approach
in~\cite{VGC} to our situation. This is possible
since the Gr\"obner basis in (\ref{redGB}) is reverse lexicographic:
first deform the lowest variable to zero,
then  deform the second lowest  variable to zero,
then the third lowest variable, and so on.

\subsection{SAGBI homotopy}\label{SAGBIhomotopy}

Let $X = (x_{ij})$ be an $(m+p) \times p$-matrix of indeterminates.
We identify the coordinate ring of the Grassmannian with the 
${\bf C}$-algebra generated by the $p \times p$-minors of $X$. 
Call this algebra $R$ and
write $[\alpha]( x_{ij} )$ for the minor indexed by $\alpha$.
Reinterpreting classical results in~\cite[\S XIV.9]{Hodge_Pedoe}, it was
shown in~\cite[Theorem 3.2.9]{Sturmfels_invariant} 
that these generators form a {\it SAGBI basis} with respect to
the lexicographic term order induced from $\,x_{11} > x_{12} >
\cdots > x_{1p} > x_{21} > \cdots > x_{m+p,p}$. This means that
the {\it initial algebra} ${\bf C}[\,in_>(f)\,:\, f \in R \,]\,$
is generated by the main diagonal terms
of the $p \times p$-minors,
\begin{equation}\label{inR}
 {\it in}_> \bigl([\alpha](x_{ij}) \bigr) \quad = \quad
x_{\alpha_1,1}  \,x_{\alpha_2,2} \,x_{\alpha_3,3} \cdots
\cdots x_{\alpha_m,m} .
\end{equation}
The resulting flat deformation can be realized by 
replacing $x_{ij}$ with $\,x_{ij}  t^{(i-1)(p-j)} \,$
for $t \rightarrow 0$ in the matrix $X$. If we expand
$\,[\alpha]( x_{ij} t^{(i-1)(p-j)}) \,$ as a polynomial
in $t$, then the lowest term equals $t^{w_a}$ times
the main diagonal monomial (\ref{inR}), where
$$ 
w_\alpha \quad := \quad 
\sum_{j=1}^p \, (\alpha_j - 1)(p - j) . 
$$
In what follows we multiply that polynomial by
$\,t^{-w_\alpha}$. For any $t \in {\bf C}$ consider the algebra
$$ 
R_t \quad := \quad {\bf C} \left[
\,\, t^{-w_\alpha} \cdot [\alpha]( x_{ij}  t^{(i-1)(p-j) } ) \,\,:\,\
\alpha \in  {\textstyle \binom{[m+p]}{p}}   \,\right].
$$
Then $R_1$ is the coordinate ring of the Grassmannian, and
$\,R_0\,$ is the algebra generated by the monomials (\ref{inR}).
This is a flat deformation of ${\bf C}$-algebras;
see~\cite{CHV} and~\cite[\S 11]{Sturmfels_GPCP}.

The {\it SAGBI homotopy} is the following system of $mp$ equations:
\begin{equation}\label{SAGBIh}
 \sum_{\alpha \in \binom{[m+p]}{p}} C^{i}_\alpha \cdot
t^{-w_\alpha} \cdot [\alpha]( x_{ij}  t^{(i-1)(p-j) } ) 
 \quad = \quad 0 \quad \qquad (i=1,\ldots,mp). 
\end{equation}
We reduce the number of variables to $mp$ by introducing
local coordinates as follows:
$\,x_{ii} = 1 \,$ for $i=1,\ldots,p$ and
$\,x_{ij} = 0 \,$ for $i<j$ or $i > m+j $.
Our original problem is to solve the system (\ref{SAGBIh}) for $t=1$.

The flatness of the family of algebras $R_t$ guarantees that the
system (\ref{SAGBIh}) has the same finite number of complex solutions
(counting multiplicities) for almost every $t \in {\bf C}$.
For $t=0$ we get a system of linear  equations in $R_0$:
\begin{equation}\label{R0eqs}
 \sum_{\alpha \in \binom{[m+p]}{p}} C^{i}_\alpha \cdot
x_{\alpha_1,1}  x_{\alpha_2,2} x_{\alpha_3,3} \cdots x_{\alpha_p,p}
\quad \qquad (i=1,\ldots,mp). 
\end{equation}

In order to solve these equations we apply the
symbolic-numeric algorithm in~\cite{Huber_Sturmfels}, while
taking advantage of the following combinatorial structures
described in~\cite[Remark 11.11]{Sturmfels_GPCP}. 
The common Newton polytope of the equations (\ref{R0eqs})
equals the {\it order polytope} of the product of an
$m$-chain and a $p$-chain (Sturmfels, Remark 11.11).
We have the following combinatorial result.

\begin{prop}
The following five numbers coincide:
\begin{itemize}
\item the right hand side of  (\ref{grassdeg}),
\item the number of linear extensions of the product of a
$m$-chain and a $p$-chain,
\item the number of maximal chains in Young's poset,
\item the normalized volume of the order polytope, and
\item the number of roots in $({\bf C}^*)^{mp}$ of
a generic system (\ref{R0eqs}).
\end{itemize}
\end{prop}

The equality of (4) and (5) is a special case of Kouchnirenko's 
Theorem~\cite{Kouchnirenko}, as all the equations have the same Newton
polytope.
The order polytope has a distinguished unimodular regular
triangulation with simplices indexed by the chains in Young's poset.
This regular triangulation is induced by the system of weights
given in (\ref{weights}). We may use these weights to define a numerical
homotopy for finding all isolated solutions of (\ref{R0eqs}).
Once this is accomplished, we trace these roots from $t=0$ to $t=1$
in the homotopy~(\ref{SAGBIh}).

\section{Special Schubert conditions}\label{pieri-deformations}


We first  describe a purely combinatorial
method for computing the number $d$ of solution planes.
Instead of the algebraic relation (1.2) we shall make use
of Young's poset which was introduced in Section 2.2.
A cover $[\alpha]\lessdot [\beta]$ in Young's poset
determines a unique index $j=j(\alpha,\beta)$ for which
$$  \alpha_j +1 \ =\ \beta_j \quad\mbox{and}
\quad\alpha_i\ =\ \beta_i \ \ \ \mbox{for} \ \ \ i\neq j. $$
A chain $ \,[\alpha^0]\lessdot [\alpha^1]\lessdot \cdots\lessdot[\alpha^l]\,$ 
in Young's poset is {\em increasing at $i$} if either $i=1$, or else
$i>1$ and $j(\alpha^{i-2},\alpha^{i-1})\leq j(\alpha^{i-1},\alpha^i)$.
For instance, $[12]\lessdot [13]\lessdot [14]\lessdot [24]$
is increasing at 1 and 2, but decreasing at 3.

Given positive integers $r_0,\ldots,r_a$, the {\em Pieri tree} 
${\mathcal T}(r_0,\ldots,r_a)$ consists of all chains 
of length $r_0+\cdots+r_a$ in Young's poset
which start at the bottom element $[1,2,\ldots,p]$ and 
which increase everywhere, except possibly
at $r_0+1,r_0+r_1+1,\ldots,r_0+\cdots+r_{a-1}+1$.
Here we include all initial segments of such chains and
we order the chains by inclusion.
Label a node in the Pieri tree by the
endpoint of the chain which that node represents.
Then the sequence of labels from the root to that
node is the chain which that node represents. 
For example, here is ${\mathcal T}(2,2)$ when $m=5$, $p=2$:
$$
\epsfxsize=1.05in\epsfbox{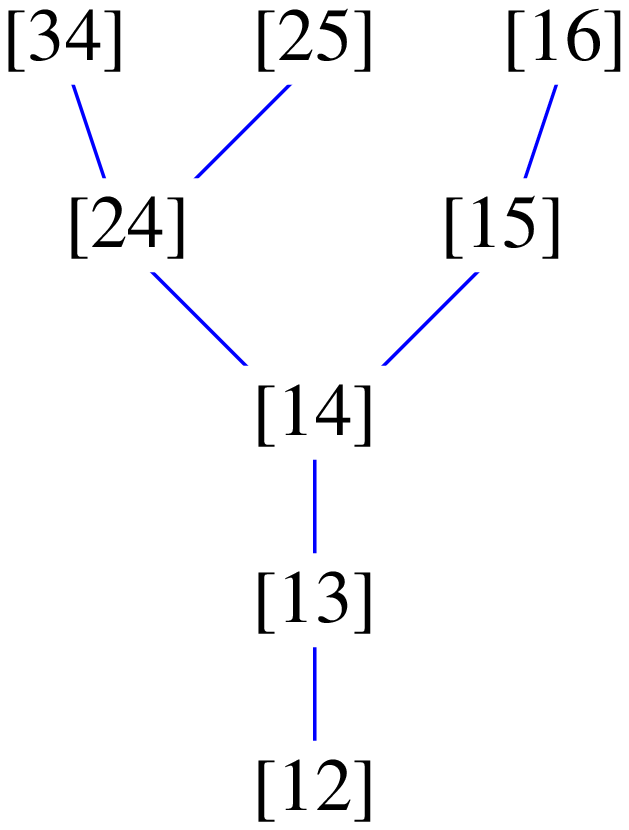} 
$$

To compute the number $d$, partition the integer sequence
 $k_1,\ldots,k_n$ into three parts
$r_0,\ldots,r_a$, $r'_0,\ldots,r'_{a'}$, and $q$.
Then $d$ is the number of pairs $(R,S)$ where $R$ is a 
leaf of ${\mathcal T}(r_0,\ldots,r_a)$, $S$ is a leaf of 
${\mathcal T}(r'_0,\ldots,r'_{a'})$, and 
the endpoints $[\alpha]$ of $R$ and $[\alpha']$ of $S$ satisfy
{\it Pieri's condition}:
\begin{equation}\label{eq:pieri}
\alpha'_1\leq m+p+1-\alpha_p <\alpha'_2\leq 
\cdots<\alpha'_p\leq m+p+1-\alpha_1
\end{equation}
Call this set of pairs 
$\mbox{\em Sols} = \mbox{\em Sols}(r_0,\ldots,r_a;r'_0,\ldots,r'_{a'})$.

For instance, $d=6$ for the sequence $2,2,2,2,2$ with $(m,p)=(5,2)$:
Of the 9 pairs $(R,S)$ of leaves of ${\mathcal T}(2,2)$,
only the following 6 satisfy~(\ref{eq:pieri}): 
(Here we  represent a leaf by its label.)
\begin{equation}\label{Sols}
\left\{
([25],[34]),\ ([34],[25]),\ ([25],[25]),\ ([25],[16]),\ ([16],[25]),
\ ([16],[16]) \right\}. 
\end{equation}

This combinatorial rule for the number $d$ gives the same answer as
the algebraic rule (1.2)
because the Pieri tree and Pieri's condition~(\ref{eq:pieri}) encode the
structure of the cohomology ring ${\mathcal A}_{m,p}$ with 
respect to its Schur basis
$\{S_\lambda \mid p\geq \lambda_1\geq \cdots\geq
\lambda_{p+1}=0\}$; see~\cite[\S I]{Macdonald_symmetric} 
or~\cite[\S 9.4]{Fulton_tableaux}.  
Specifically,
$$ 
\prod_{i=0}^a h_{r_i} \ =\ 
\sum_{\lambda_1+\cdots +\lambda_p=r_0+\cdots+r_a}
K_{\lambda,(r_0,\ldots,r_a)}\cdot S_\lambda,
$$
where $K_{\lambda,(r_0,\ldots,r_a)}$ is the number of leaves in 
${\mathcal T}(r_0,\ldots,r_a)$ with label
$$
[\alpha(\lambda)]\quad  :=\quad 
[\lambda_p+1,\lambda_{p-1}+2,\ldots,\lambda_1+p].
$$
The numbers $K_{\lambda,(r_0,\ldots,r_a)}$ are called
{\it Kostka numbers}. In ${\mathcal A}_{m,p}$ we calculate
\begin{eqnarray*}
\prod_{i=1}^n h_{k_i} &=&
\left(\prod_{i=0}^a h_{r_i}\right)\cdot
\left(\prod_{j=0}^{a'} h_{r'_j}\right)\cdot h_q\\
&=&
\sum_{\lambda,\mu}K_{\lambda,(r_0,\ldots,r_a)} 
K_{\mu,(r'_0,\ldots,r'_{a'})}\cdot S_\lambda\cdot S_\mu\cdot h_q.
\end{eqnarray*}
We evaluate this expression with Pieri's 
formula (Proposition~\ref{prop:triple} below):
If $\lambda_1+\cdots+\lambda_p+\mu_1+\cdots+\mu_p+q=mp$,
then $S_\lambda\cdot S_\mu\cdot h_q$ is either $(h_m)^p$  or 0 depending
upon whether or not
$$
\lambda_p\leq m-\mu_1 \leq \lambda_2\leq\cdots 
\leq \lambda_1 \leq m-\mu_p. \eqno (\ref{eq:pieri}')
$$
The condition~(\ref{eq:pieri}$'$)
is equivalent to~(\ref{eq:pieri}) under the transformation 
$\,\lambda\leftrightarrow[\alpha(\lambda)]$.

These methods correctly enumerate the $p$-planes which meet
each $K_1,\ldots,K_n$ nontrivially
because, under the isomorphism between
${\mathcal A}_{m,p}$ and the cohomology ring of 
$\mbox{\em Grass}(p,m+p)$, the indeterminate
$h_{k_i}$ corresponds to the cohomology class Poincar\'e dual to 
$\Omega_{K_i}$, the set of $p$-planes which meet $K_i$ nontrivially.
Moreover,  $(h_m)^p$ represents the class dual to a point.
The Pieri tree models certain intrinsic deformations (described in 
\S 3.3 and \S \ref{explicit_pieri}) of the Grassmannian
which establish this isomorphism, and which we shall use
for computing the $p$-planes which meet each $K_1,\ldots,K_n$ nontrivially. 

\subsection{Basics on Schubert varieties}
For vectors $f_1,\ldots,f_j$ in ${\bf C}^{m+p}$, let 
$\langle f_1,\ldots,f_j\rangle$ be their linear span.
Fix the columns $e_1,\ldots,e_{m+p}$ of the identity matrix as a 
standard basis for ${\bf C}^{m+p}$.
For $\alpha\in\binom{[m+p]}{p}$, define 
$\alpha^\vee\in\binom{[m+p]}{p}$ by 
$\alpha^\vee_j := m+p+1-\alpha_{p+1-j}$.

A sequence $\alpha\in\binom{[m+p]}{p}$ determines a {\em Schubert variety}
$$
\Omega_\alpha\quad:=\quad \{X\in\mbox{\em Grass}(p,m+p)\, \, \mid \,\,
\dim X\cap\langle e_1,\ldots,e_{\alpha^\vee_j}\rangle \geq j \,\,\,
\hbox{for} \,\, 1\leq j\leq p\}.
$$
This variety has complex codimension
$|\alpha|:=\alpha_1-1+\alpha_2-2+\cdots+\alpha_p-p$. 
Similarly define 
$$ 
\Omega'_\alpha\quad:=\quad 
\{X \in\mbox{\em Grass}(p,m+p)\, \, \mid \,\,
 \dim X \cap\langle e_{\alpha_j},\ldots,e_{m+p}\rangle \geq p+1-j \,\,\,\,
\hbox{for} \,\,  1\leq j\leq p\}.
$$
For a linear subspace $N$ of ${\bf C}^{m+p}$ of dimension $m+1-q$, define
the {\em special Schubert variety}
$$
\Omega_N \quad := \quad  \{X\in\mbox{\em Grass}(p,m+p) \,\, \mid \,\,
\dim X \cap N \geq 1\}.
$$
This  has codimension $q$.
If $N=\langle e_1,\ldots,e_{m+1-q}\rangle$, then 
$\Omega_N = \Omega_{[1,\ldots,p-1,p+q]}$.

The special Schubert variety $\Omega_N$ is cut out
by the system of $\binom{m+p}{q-1}$ polynomial equations:
\begin{equation}\label{maxminors}
X\ \in \ \Omega_N \quad \Longleftrightarrow \quad
\mbox{all maximal minors of }\left[ X\mid N\right]
\, \mbox{are zero,}
\end{equation}
where $X\in \mbox{\em Grass}(p,m+p)$ is represented by a 
$(m+p)\times p$-matrix.
The Laplace expansion of these  equations in terms of the
Pl\"ucker coordinates of $X$ define $\Omega_N$ as a subscheme of 
$\mbox{\em Grass}(p,m+p)$.
These equations are redundant: select
$ m\! + \! 1 \! - \! q $ rows of
$\left[ X\mid N\right]$  such that the corresponding
maximal minor of $N$ is invertible.
Consider the set of maximal minors of
$\left[ X\mid N\right]$ which cover all the rows selected.
This gives $\binom{p}{q-1}$ polynomial equations which 
generate the same ideal as all  $\binom{m+p}{q-1}$
minors of $\left[ X\mid N\right]$.
For a purely set-theoretic (but not scheme-theoretic)
representation of $\Omega_N$ a further substantial reduction
in the number of equations is possible using the
results of~\cite{Bruns_Schwanzl}.

An intersection $Y\cap Z$ of subvarieties is {\em generically transverse}
if every component of  $Y\cap Z$ has an open subset along which $Y$ and
$Z$ meet transversally.
In this case the following identity in the cohomology ring
holds:
$$ [Y \cap Z]\quad = \quad [Y]\cdot[Z], $$
where $[W]$ denotes the cycle class of a subvariety $W$.
By Kleiman's Transversality Theorem~\cite{Kleiman},
subvarieties of ${\it Grass}(p,m+p)$ in general position meet generically
transversally.
Transversality and generic transversality coincide when $Y\cap Z$ 
is finite.

\begin{prop}[Hodge and Pedoe, 1952, Theorem III in 
\S XIV.4]\label{prop:triple}  \hfill \break
Let $\alpha,\,\alpha'\in\binom{[m+p]}{p}$ with $|\alpha|+|\alpha'|+q=mp$
and let $N$ be a linear subspace of ${\bf C}^{m+p}$ 
with dimension $m+1-q$ none of whose Pl\"ucker coordinates vanish.
Then the intersection
\begin{equation}\label{triple_int}
\Omega_\alpha\,\cap\,\Omega'_{\alpha'} \, \cap \, \Omega_N
\end{equation}
either is transverse consisting of a single $p$-plane or is empty,
depending upon whether or not~(\ref{eq:pieri}) holds.
\end{prop}

\noindent{\sc Proof and Algorithm: } 
The intersection $\Omega_\alpha \cap\Omega'_{\alpha'}$
is nonempty if and only if $\alpha'_j\leq \alpha^\vee_j$
for $j = 1,\ldots,p$. These are the weak inequalities
in~(\ref{eq:pieri}). We shall assume that they hold  in what follows.
The $p$-planes in $\Omega_\alpha \cap\Omega'_{\alpha'}$ are represented
by $(m+p)\times p$-matrices $X=(x_{ij})$ such that
\begin{equation}\label{eq:alpha}
x_{i,j}\ =\ 0\qquad\mbox{for}\qquad 
i<\alpha'_j \quad \mbox{or} \quad
\alpha^\vee_j < i .
\end{equation}
Consider the nonzero coordinate subspaces $\,C_j
:=\langle e_{\alpha'_j},\ldots,e_{\alpha^\vee_j}\rangle$, set
$\,C := C_1 + \cdots + C_p $, and note that
$\, p+q = \sum_j \dim (C_j) \geq dim(C)$.
{}From~(\ref{eq:alpha}) we see that  $\,X 
\in \Omega_\alpha\cap \Omega'_{\alpha'}\,$ implies
$\, X\subseteq C \,$ and hence $\, N \cap X \subseteq N \cap C$.
Therefore the triple intersection~(\ref{triple_int})
is nonempty only if the following equivalent conditions hold:
$$
 dim(C\cap N) \geq 1 \,\, \iff\,\,
  dim(C) = p+q \,\,
 \iff \,
 \mbox{the sum}  \,\,\,C = C_1 + \cdots + C_p \,\,\, \mbox{is direct} $$
$$
\iff \,
\alpha^\vee_j < \alpha'_{j+1} \,\,\mbox{for} \,\, j=1,\ldots,p-1 \, \iff \,
\mbox{(3.1) holds.} $$
In this case we determine$\,C \cap N \,$ by computing vectors
$\,g_j\in C_j \,$  such that
$\,C \cap N = \langle g_1\oplus g_2\oplus\cdots\oplus g_p \rangle $.
(This computation is the ``algorithm'' part in this proof.)
The desired $p$-plane $X$ satisfies $\, C \cap N = X \cap N$, and,
in view of~(\ref{eq:alpha}), this implies $X = \langle g_1,\ldots,g_p\rangle$.
Transversality of~(\ref{triple_int})
is verified in local coordinates for
$\Omega_\alpha\cap \Omega'_{\alpha'}$ by considering $p+q-1$ independent
linear forms which vanish on $N$.\qed

For $\alpha\in \binom{[m+p]}{m}$ define 
$\lambda(\alpha)$ 
by $\lambda(\alpha)_j:= \alpha_j-j$ for $1\leq j\leq p$.
Then $S_{\lambda(\alpha)}$ represents the cycle class 
of $\Omega_\alpha$ (equivalently, of $\Omega'_\alpha$).
If $\dim N=m+1-q$, then  $h_q$ is the cycle class of $\Omega_N$.
Suppose $|\alpha|+|\alpha'|+q=mp$.
Then Proposition~\ref{prop:triple} implies the following identity
in ${\mathcal A}_{m,p} $:
$$ S_{\lambda(\alpha)}\cdot S_{\lambda(\alpha')}\cdot h_q\ =\ 
\left \{\begin{array}{cl}(h_m)^p&\quad\mbox{if } (\ref{eq:pieri})
\mbox{ holds} \\
0&\quad\mbox{otherwise}\end{array}\right.. $$
This identity implies (via Poincar\'e duality) that 
$$ S_{\lambda(\alpha)}\cdot h_q \quad = \quad \sum S_{\lambda(\beta)},
\label{pieri_product} 
$$
the sum over all $\beta$ with $|\beta|=|\alpha|+q$ for which
$\alpha,\beta^\vee$ satisfy Pieri's condition~(\ref{eq:pieri}).
Call this set $\alpha*q$, which is 
also the set of endpoints of increasing chains of length $q$ in Young's
poset that begin at $\alpha$.

This last form has geometric content.
In~\cite{Sottile_explicit_pieri}, explicit deformations were given that
transform the irreducible intersection 
$\Omega_\alpha \cap \Omega_N$ into the cycle
$\sum_{\beta\in\alpha*q} \Omega_\beta$.
Moreover, the branching of the components of the cycles in these
deformations reflects the branching among these increasing chains above
$\alpha$.
This process may be iterated 
to transform  an intersection of several special Schubert varieties into a
sum of triple intersections of the form~(\ref{triple_int}),
indexed by pairs $(R,S)\in\mbox{\em Sols}$.
{}From this sum, we obtain a set of start solutions indexed by {\em Sols}.
Also, every intermediate cycle in these deformations consists of the same
number (counting multiplicities) of $p$-planes.
The Pieri homotopy begins with one of the start solutions and uses
numerical path continuation to trace the
sequence of curves defined by these deformations which connect that start
solution to a solution of the original problem.

\subsection{Pieri homotopy algorithm}\label{pieri_alg}
Given linear subspaces $K_1,\ldots,K_n$ in general position with 
$\dim K_i=m+1-k_i$ and  
$k_1+\ldots+k_n=mp$, first partition  $K_1,\ldots,K_n$ into three
lists:
$$
L_0,\ldots,L_a,\qquad 
L'_0,\ldots,L'_{a'},\qquad 
N
$$
where $\dim L_i=m+1-r_i$, $\dim L'_i=m+1-r'_i$, and 
$\dim N=m+1-q$.
Construct the Pieri trees ${\mathcal T}(r_0,\ldots,r_a)$ and 
${\mathcal T}(r'_0,\ldots,r'_{a'})$, and form the set  {\it Sols}.
Change coordinates so that $L_0=\langle e_1,\ldots,e_{m+1-r_0}\rangle$
and $L'_0 =\langle e_{p+r'_0},\ldots,e_{m+p}\rangle$.
Set $\tau := \max \{r_1+\cdots+r_a,\ r'_1+\cdots+r'_{a'}\}$.

Given a chain $R$ in the Pieri tree and a positive integer $k$, let
$R(k)$ be the $k$th 
element in that chain, or, if $k$ exceeds the length of $R$, then let
$R(k)$ be the endpoint of $R$.
For each $(R,S)\in{\it Sols}$ and $k$ from $\tau$ to $0$
we shall construct (in Definition~\ref{def:ZRkt} below) one-parameter
families  $Z_{R,k}(t)$ and $Z'_{S,k}(t)$ of pure-dimensional
subvarieties of  ${\it Grass}(p,m+p)$ with the following properties:
\begin{enumerate}
\item $Z_{R,k}(t)\subset \Omega_{R(r_0+k)}$ and 
$Z'_{S,k}(t)\subset \Omega'_{S(r'_0+k)}$.
\item For $t=0$ or $1$ and each $k$,
$Z_{R,k}(t)\cap Z'_{S,k}(t)\cap \Omega_N$
is transverse and 0-dimensional.
\item $Z_{R,\tau}(t) = \Omega_{R(r_0+\tau)}$ and 
$Z'_{S,\tau}(t) = \Omega'_{S(r_0+\tau)}$ .
\item 
$Z_{R,k+1}(1)$ is a component of $Z_{R,k}(0)$.
Likewise, $Z'_{S,k+1}(1)$ is a component of $Z'_{S,k}(0)$.
\item
$Z_{R,0}(1) = \Omega_{L_0}\cap\cdots\cap\Omega_{L_a}$ and 
$Z'_{S,0}(1) = \Omega_{L'_0}\cap\cdots\cap\Omega_{L'_{a'}}$.
\end{enumerate}

Property 4 is a consequence of 
Proposition~\ref{prop:pieri}, the others follow from the assumption of
genericity and the definition (Definition~\ref{def:ZRkt}) of the
families $Z_{R,k}(t)$ and $Z'_{S,k}(t)$. 

By 2, the family $W_{(R,S),k}(t)$ over ${\bf C}$ whose
fibre at general $t$  (including $t=0$ and $t=1$) is 
$$
W_{(R,S),k}(t) \quad:=\quad Z_{R,k}(t)\cap Z'_{S,k}(t)\cap\Omega_N
$$
consists of a finite number of curves.
In fact, for general $t$  (including $t=0$ and $t=1$), 
$W_{(R,S),k}(t)$ has the following general form
(see Definition~\ref{def:ZRkt} for the precise form):
$$
W_{(R,S),k}(t)\ =\ 
\Omega_\alpha\cap\Omega'_{\alpha'}\cap
\Omega_{M_1}\cap\cdots\cap\Omega_{M_s},
$$
where $M_1,\ldots, M_s$ are linear subspaces with $M_s=N$, which depend
upon $R,S,k$, 
the subspaces $L_1,\ldots,L_a$, $L'_1,\ldots,L'_{a'}$, and at most two of
the $M_i$ depend upon $t$.
Also, $\alpha$ and $\alpha'$ depend upon $R,S$, and $k$
with the typical case being $\alpha= R(r_0+k)$ and  $\alpha'= S(r'_0+k)$.

The numerical homotopy defined by the curves $W_{(R,S),k}(t)$ may be
expressed in a parameterization $X=(x_{i,j})$ of an open subset of
$\Omega_\alpha\cap\Omega'_{\alpha'}$:
\begin{equation}\label{loccord}
x_{i,j}\ =\ 0\qquad\mbox{if}\qquad i<\alpha'_j \quad\mbox{or}\quad
\alpha^\vee_j<i\quad\qquad\mbox{and}\quad\qquad x_{\delta_j,\,j}\ =\ 1,
\end{equation}
where  $\delta:=S(r'_0+\tau)$.
The equations for  $W_{(R,S),k}(t)$ are then
$$
\mbox{maximal minors }\left[ X\mid M_i\right]\ =\ 0 \qquad
i=1,\ldots,s.
$$

The curves of $W_{(R,S),k}(t)$ define the sequences of homotopies in the
Pieri homotopy algorithm as follows:
For $(R,S)\in {\it Sols}$, let $X_{(R,S),\tau}$ be the (unique by
Proposition~\ref{prop:triple}) 
$p$-plane in
$\Omega_{R(r_0+\tau)}\cap\Omega'_{S(r'_0+\tau)}\cap\Omega_N = 
W_{(R,S),\tau}(1)$. 
By 3 and 4, $X_{(R,S),\tau}\in W_{(R,S),\tau-1}(0)$ and hence lies on a unique
curve in $W_{(R,S),\tau-1}(t)$.
Use numerical path continuation to trace this curve from $t=0$ to $t=1$
to obtain $X_{(R,S),\tau-1}$, which is 
a point of $W_{(R,S),\tau-2}(0)$, by 4.
Then $X_{(R,S),\tau-1}$ lies on a unique curve in $W_{(R,S),\tau-2}(t)$,
which we trace to find $X_{(R,S),\tau-2}\in W_{(R,S),\tau-2}(1)$.
Continuing this process, after tracing $\tau$ curves, we obtain
$X_{(R,S),0}\in W_{(R,S),0}(1)$, which is a solution to the original system,
by 5.  We show shall prove in Theorem~\ref{allsols} that 
$\{X_{(R,S),0}\mid (R,S)\in {\it Sols}\}$ consists of all the solutions to the
original system.

\subsection{Definition of the moving cycles $Z_{R,k}(t)$}

The cycle $Z_{R,k}(t)$ will depend upon the choice of a general
upper triangular $(m+p)\times (m+p)$-matrix  $F$ with 1's on its
anti-diagonal, 
$$
\left(\begin{array}{ccc}*&*&1\\
{}*&\hspace{1pt}\raisebox{0pt}{.}\raisebox{2.2pt}{.}%
\raisebox{4.4pt}{.}
\\1&&0\end{array}\right),
$$
the $k$th link in the chain $R$, and the data $L_1,\ldots,L_a$.
The key ingredient of this definition of $Z_{R,k}(t)$ is the
construction of a one-parameter family of linear subspaces
$\Lambda_i(t)$ in Definition~\ref{def:Lambda}, which depends upon  $F$.
The matrix $F$ is fixed throughout the algorithm, its purpose is that 
$\langle e_{m+p-j},\ldots,e_{m+p}\rangle$ equals the span of the first 
$j$ columns of $F$, and these columns are in general position with
$e_1,\ldots,e_{m+p}$.
The subtle linear degeneracies of $\Lambda_i(t)$ as $t\rightarrow 0$ are at
the heart of this homotopy algorithm, as well as the explicit proof of
Pieri's formula~\cite[Theorem~3.6]{Sottile_explicit_pieri}, which we state
below (Proposition~\ref{prop:pieri}).

\begin{defn}\label{def:ZRkt}\mbox{\ }

\begin{enumerate}
\item If $r_1+\cdots+r_a\leq k$, then set \ 
$Z_{R,k}(t)\ :=\ \Omega_{R(r_0+\tau)}$.
\item Otherwise, define $c$ by 
$r_1+\cdots+r_{c-1}\leq k<r_1+\cdots+r_c$, and 
set $i:= k-r_1-\cdots-r_{c-1}$,
$\alpha:= R(r_0+r_1+\cdots+r_{c-1})$, and $\beta:= R(r_0+k)$.
\begin{enumerate}
\item If $i>0$ and $\beta_p>\alpha_p$, then 
$\beta+(0,\ldots,0,r_c-i) = R(r_0+\cdots+r_c)$, and we set 
$$
Z_{R,k}(t)\ :=\ \Omega_{R(r_0+\cdots+r_c)}
\cap \Omega_{L_{c+1}}\cap\cdots\cap\Omega_{L_a}.
$$
\item Otherwise, let $\Lambda_i(t)$ be the 1-parameter family of linear
subspaces given by Definition~\ref{def:Lambda}, where we let 
$L:= L_c$ and $r:= r_c$. 

If $i=0$, then $\Lambda_0(1)=L_c$, $\alpha=\beta=R(r_0+k)$, and we set
$$
Z_{R,k}(t)\ :=\ \Omega_{R(r_0+k)}\cap
\Omega_{\Lambda_0(t)}
\cap \Omega_{L_{c+1}}\cap\cdots\cap\Omega_{L_a}.
$$
If $i>0$, let $j$ be maximal such that $\beta_j>\alpha_j$.
Then $j<p$ as $j=p$ is case 2(a).
Set
$$
Z_{R,k}(t)\ :=\ \Omega_{R(r_0+k)}\cap
\Omega_{\Lambda_i(t)\cap\langle e_1,\ldots,e_{\beta^\vee_{p+1-j}}\rangle}
\cap \Omega_{L_{c+1}}\cap\cdots\cap\Omega_{L_a}.
$$
\end{enumerate}
\end{enumerate}
We define $Z'_{S,k}(t)$ similarly, but with the matrix $F$  replaced by a
lower triangular matrix with 1's on its diagonal, and 
$\langle e_1,\ldots,e_{\beta^\vee_{p+1-j}}\rangle$ replaced by
$\langle e_{\beta_j},\ldots,e_{m+p}\rangle$.
\end{defn}

\begin{defn}\label{def:Lambda}
Let $F$ be an upper triangular matrix with 1's on the anti-diagonal and 
$L$ be a general $(m+1-r)$-plane, represented as a  
$(m+p)\times(m+1-r)$-matrix with columns $l_1,\ldots,l_{m+1-r}$.
Construct a $(m+p)\times p$-matrix $U=(u_1,\ldots,u_p)$ as follows:
Reverse the last $m+p-\alpha_p$ columns of $F$, then remove the columns
indexed by $\alpha^\vee_1,\ldots,\alpha^\vee_p$.

For each $0\leq i< r$, define a one-parameter family of 
$(m+p)\times(m+1-r)$-matrices $\Lambda_i(t)$ for $t\in{\bf C}$
as follows:
\begin{enumerate}
\item
If $i=0$, then the $b$th column of $\Lambda_0(t)$ is
$t\cdot l_b + (1-t)\cdot u_b$.

\item
For $0<i< r$, the $b$th column of $\Lambda_i(t)$
is 
$$
\begin{array}{lcl}
t\cdot u_{b+i-1} + (1-t)\cdot u_{b+i}     &\ & b+i-1<\alpha_p-p\\
t\cdot u_{b+i-1} + (1-t)\cdot u_{p+1+i-r} &&   b+i-1=\alpha_p-p\\
       u_{b+i-1}		     &&   b+i-1>\alpha_p-p
\end{array}
$$
\end{enumerate}

\end{defn}

\subsection{An example}
We give an example illustrating these definitions and the Pieri
homotopy algorithm. 
Let $L_0,L_1,L'_0,L'_1$, and $N$ be general 4-planes in ${\bf C}^7$.
We give a sequence of homotopies $W_{(R,S),k}(t)$ for $k=2,1,0$
for finding one of the six $2$-planes which meet each of the five given 
4-planes nontrivially.

Here, $(m,p)=(5,2)$ and   $k_1=\cdots=k_5=2$ so that  $\tau=2$.
Construct the set {\em Sols} as in (\ref{Sols}).
Let $(R,S)\in\mbox{\em Sols}$ be the following two sequences:
$$
R\ :=\ [12]\lessdot[13]\lessdot[14]\lessdot[24]\lessdot[25],\qquad
S\ :=\ [12]\lessdot[13]\lessdot[14]\lessdot[15]\lessdot[16].
$$
Let $e_1,\ldots,e_7$ be the columns of a $7\times 7$-identity matrix,
a basis for ${\bf C}^7$.
Suppose that  $L_0=\langle e_4,e_5,e_6,e_7\rangle$ and 
$L'_0=\langle e_1,e_2,e_3,e_4\rangle$
and represent $L_1$ and $L'_1$ as $7\times 4$-matrices.
Then 
$\Omega_{[14]}\cap\Omega_L\cap\Omega'_{[14]}\cap 
\Omega_{L'}\cap\Omega_N$ is the set of 2-planes which meet all
five linear subspaces nontrivially.

We first find the plane 
$X_{(R,S),2}\in \Omega_{[25]}\cap\Omega'_{[16]}\cap\Omega_N$,
using the algorithm in the proof of Proposition~\ref{prop:triple}.
Suppose that $N$ has the form
$$
\left[ \begin{array}{c} n\\\hline  I\end{array}\right],
$$
where $I$ is the $4\times 4$ identity matrix, and $n$ is a $3\times
4$-matrix.
In this case, $C_1 = \langle e_1,e_2,e_3\rangle$ and 
$C_2= \langle e_6\rangle$, hence $C = \langle e_1,e_2,e_3,e_6\rangle$.
Thus the intersection $C\cap N$ is generated by the third
column of $N$, and so $X_{(R,S),2}$ is represented by the matrix:
$$
X_{(R,S),2}\quad=\quad
\left[\begin{array}{cc}
         n_{13} & 0 \\
         n_{23} & 0 \\
         n_{33} & 0 \\
           0    & 0 \\
           0    & 0 \\
           0    & 1 \\
           0    & 0 \\
\end{array}\right].
$$
Following Definition~\ref{def:Lambda}, we have:
\begin{eqnarray*}
\Lambda_0(t) &=&\langle tl_1 + (1-t)u_1,  tl_2 + (1-t)u_2,  
tl_3 + (1-t)u_3,  tl_4+(1-t)u_4\rangle,\\
\Lambda_1(t) &=&\langle tu_1+(1-t)u_2,tu_2+(1-t)u_5,u_3,u_4\rangle.
\end{eqnarray*}
$\Lambda'_i(t)$ is defined similarly.

We describe the families $W_{(R,S),k}(t)$ for $k=2,1,0$ in local
coordinates for  $\Omega_{[14]}\cap \Omega'_{[14]}$ 
determined by the sequence $[16]$: 
$$
 X \quad = \quad
\left[\begin{array}{cc}
             1  &    0   \\
         x_{21} &    0   \\
         x_{31} &    0   \\
         x_{41} & x_{42} \\
           0    & x_{52} \\
           0    &    1   \\
           0    & x_{72} \\
\end{array}
\right] \ .
$$
The family $W_{(R,S),2}(t)$ is the constant family
$\{X_{(R,S),2}\} = \Omega_{25}\cap\Omega'_{16}\cap \Omega_N$.
Assuming that $n_{13}$, which is the $[1457]$th Pl\"ucker coordinate of 
$N$,  is non-zero, then $X_{(R,S),2}$ may be expressed in these local
coordinates: 
$$
 X_{(R,S),2} \quad = \quad
\left[\begin{array}{cc}
         1       &  0  \\
   n_{23}/n_{13} &  0  \\
   n_{33}/n_{13} &  0  \\
         0       &  0  \\
         0       &  0  \\
         0       &  1  \\
         0       &  0  \\
\end{array}
\right] \ .
$$

When $k=1$, first consider the definition of $Z_{R,1}(t)$.
Here we are in case 2(b) with $\beta=[24]$ and $i>0$, so that 
$\beta^\vee=[46]$. 
Since $\Lambda_1(t)\subset \langle e_1,\ldots,e_6\rangle$, 
we have 
$Z_{R,1}(t)\ =\  \Omega_{[24]}\cap  \Omega_{\Lambda_1(t)}$.
For the definition of $Z'_{S,1}(t)$, we are in case 2(a), so that
$Z'_{S,1}(t)\ =\  \Omega'_{[16]}$.
Hence 
$$
W_{(R,S),1}(t)\ =\ \Omega_{[24]}\cap\Omega_{\Lambda_1(t)}
\cap\Omega'_{[16]}\cap \Omega_N.
$$
This has 3 linear equations
$x_{42}  = x_{52}  = x_{72}  = 0$,
which describe $\Omega_{[24]}\cap\Omega'_{[16]}$,
and 7 non-trivial equations, the vanishing of the maximal minors
of $[ X\mid \Lambda_1(t)]$ and $[ X\mid N]$, 
which describe $\Omega_{\Lambda_1(t)}\cap \Omega_N$.

For $k=0$, $i=0$ and we are  in case 2(b)  for both $Z_{R,0}(t)$ and
$Z'_{S,0}(t)$ so that 
$$
W_{(R,S),0}(t)\ =\ \Omega_{[14]}\cap\Omega_{\Lambda_0(t)}
\cap\Omega'_{[14]}\cap\Omega_{\Lambda'_0(t)}\cap \Omega_N.
$$
This has 21 non-trivial equations, the vanishing of the maximal
minors of $\left[ X\mid \Lambda_0(t)\right]$, 
$\left[ X\mid \Lambda'_0(t)\right]$, and $\left[ X\mid N\right]$.

\subsection{Proof of correctness}\label{explicit_pieri}

We describe the Pieri deformations
linking the families $Z_{R,k}(t)$ for $k$ from 
$r_1\!+\!\cdots\!+\!r_{c-1}$ to $r_1+\cdots+r_c$,
which establishes Property 4 of $Z_{R,k}(t)$ in \S\ref{pieri_alg}.
We also show that  the set 
$\{X_{(R,S),0}\mid (R,S)\in {\it Sols}\}$ consists of all the solutions
to the original system. 

Consider the dynamic part of $Z_{R,k}(t)$,
namely whichever of 
$$
\Omega_{R(r_0+\cdots+r_c)},\quad
\Omega_{R(r_0+k)}\cap \Omega_{\Lambda_0(t)}, \quad\mbox{or}\quad
\Omega_{R(r_0+k)}\cap
\Omega_{\Lambda_i(t)\cap\langle e_1,\ldots,e_{\beta^\vee_{p+1-j}}\rangle},
$$
appeared in the definition of $Z_{R,k}(t)$.
We call this cycle $Y_{\alpha,\beta,L}(t)$, where
$L:= L_c$ and $\alpha=R(r_0+\cdots+r_{c-1})$, and 
$\beta=R(r_0+k)$.

For $\beta\in \alpha*i$ and $\gamma \in \alpha*(i+1)$ write
$\beta\prec_\alpha\gamma$ if $\gamma$ covers $\beta$ and 
$j(\beta,\gamma)\geq j(\alpha,\beta):=
\max\{j\mid\beta_j>\alpha_j\}$.
This partitions $\alpha*(i+1)$ into sets 
$\{\gamma\mid \beta\prec_\alpha\gamma\}$ for
$\beta\in\alpha*i$.

\begin{prop}\cite[Theorem~3.6]{Sottile_explicit_pieri}\label{prop:pieri}
\mbox{ } 

Let $\alpha,\beta,i,r,L,\Lambda_i(t)$, and $Y_{\alpha,\beta,L}(t)$ be as
above.
Then
\begin{enumerate}
\item  For all $t$, $\Omega_\alpha\cap \Omega_{\Lambda_0(t)}$ is
generically transverse.
\item $Y_{\alpha,\beta,L}(t)$ is free of
multiplicities for all $t$ and irreducible for $t\neq 0$.
\item If $i\neq r-1$, then 
$Y_{\alpha,\beta,L}(0)=\sum_{\beta\prec_\alpha\gamma}
Y_{\alpha,\gamma,L}(1)$.
\item If $\beta\in \alpha*(r-1)$, then 
$Y_{\alpha,\beta,L}(0)=\sum_{\beta\prec_\alpha\gamma}
\Omega_\gamma$.
\end{enumerate}
\end{prop}

By 3, the cycle class of 
$\sum_{\beta\in \alpha*i}Y_{\alpha,\beta,L}(t)$
is independent of $i$ and $t$, and it equals the cycle class of 
$\Omega_\alpha\cap \Omega_{\Lambda_0(1)}=\Omega_\alpha\cap \Omega_L$.
By 4, we see that the cycle classes of $\Omega_\alpha\cap \Omega_L$
and $\sum_{\beta\in \alpha*r}\Omega_\beta$ coincide, furnishing another
proof of Pieri's formula.
Property 4 of $Z_{R,k}(t)$ follows from assertion 3.

\begin{thm}\label{allsols}
When $K_1,\ldots,K_n$ are generic, the Pieri homotopy algorithm finds
all $p$-planes which meet each $K_1,\ldots,K_n$ nontrivially. That is,
$$ \{X_{(R,S),0}\mid (R,S)\in {\it Sols}\}\quad=\quad
\Omega_{K_1}\cap\cdots\cap\Omega_{K_n}. $$
\end{thm}

\noindent{\bf Proof. }
Note that for any $R,S\in{\it Sols}$, the families
$Z_{R,k}(t), Z'_{S,k}(t)$, and $W_{(R,S),k}(t)$ depend only upon the
initial segments $R(0),\ldots,R(r_0+k)$ and $S(0),\ldots,S(r'_0+k)$ of
$R$ and $S$. 

By construction, the original system
$\Omega_{K_1}\cap\cdots\cap\Omega_{K_l}$ coincides with 
$W_{(R,S),0}(1)$, for any $(R,S)\in\mbox{\em Sols}$.
We inductively construct chains $R\in {\mathcal T}(r_0,\ldots,r_a)$ and 
$S\in {\mathcal T}(r'_0,\ldots,r'_{a'})$, and $p$-planes $X_k$ for 
$0\leq k\leq \tau$ such that 
$$
X_{k+1}\ \in \ W_{(R,S),k}(0)\cap W_{(R,S),k+1}(1)
$$
and $X_k, X_{k+1}$ lie on the same curve of $W_{(R,S),k}(t)$.
Then $X_\tau$ is the start solution $X_{(R,S),\tau}$, which shows that 
$X_0\in\{X_{(R,S),0}\mid (R,S)\in {\it Sols}\}$.

First set $R(0),\ldots, R(r_0)$ to be the unique chain from 
$[1,\ldots,p]$ to $[1,\ldots,p-1,p+r_0]$, and similarly for
$S(0),\ldots,S(r'_0)$. 
Then $X_0\in W_{(R,S),0}(1)$ and hence lies on a unique curve in
$W_{(R,S),0}(t)$. 
Let $X_1$ be the point on that curve with $t=0$.
By Proposition~\ref{prop:pieri} (4), 
$$
X_1\ \in\ Y_{(R(r_0),R(r_0),L_1)}(0)\ =\ 
\sum_{\beta\in R(r_0)*1} Y_{(R(r_0),\beta,L_1)}(1).
$$
Let $R(r_0+1)$ be the index $\beta$ such that 
$X_1\in Y_{(R(r_0),\beta,L_1)}(1)$.
Define $S(r'_0+1)$ similarly.
Then $X_1\in W_{(R,S),1}(1)$.

In general, suppose that we have constructed $R(0),\ldots,R(r_0+k)$,
$S(0),\ldots,S(r'_0+k)$, and $X_k\in W_{(R,S),k}(1)$.
Then $X_k$ lies on a unique curve in $W_{(R,S),k}(t)$.
Let $X_{k+1}$ be the point on that curve at $t=0$.
Let $c$ be minimal subject to $k<r_1+\cdots+r_c$ and set
$\alpha=R(r_0+\cdots+r_{c-1})$ and 
$\beta=R(r_0+k)$.
If $k+1<r_1+\cdots+r_c$, then 
by Proposition~\ref{prop:pieri} (3), there is a unique index 
$\gamma\in \beta*1$ such that 
$X_{k+1}\in Y_{\alpha,\gamma,L_c}(1)$.
If $k+1=r_1+\cdots+r_c$ then 
by Proposition~\ref{prop:pieri} (4),
there is a unique index $\gamma\in \beta*1$
such that $X_{k+1}\in\Omega_\gamma$.
Set $R(r_0+k+1)=\gamma$ and likewise define $S(r'_0+k+1)$.
Continuing in this fashion, we construct the chains $R$ and $S$, and $X_j$
for $0\leq j\leq \tau$.

We show that $R$  increases everywhere, except possibly at 
$1,r_0+1,\ldots,r_0+\cdots+r_{a-1}+1$, and hence 
$R\in{\mathcal T}(r_0,\ldots,r_a)$.   
Similar arguments show that $S\in{\mathcal T}(r'_0,\ldots,r'_{a'})$,
which will complete the proof.
Suppose $k+1\not\in\{1,r_1+1,\ldots,r_1+\cdots+r_{a-1}+1\}$.
Let $c$ be minimal subject to $k<r_1+\cdots+r_c$ and let
$\alpha=R(r_0+\cdots+r_{c-1})$, $\beta=R(r_0+k)$, and $\gamma=R(r_0+k+1)$.
Then by Proposition~\ref{prop:pieri} (3) and (4), $\beta\prec_\alpha\gamma$.
The condition $j(\beta,\gamma)\geq j(\alpha,\beta)$ in the definition of 
$\beta\prec_\alpha\gamma$ ensures that $R$ increases at $k+1$.
\qed

\section{Homotopy continuation of overdetermined systems}

Numerical homotopy continuation is a method for finding
the isolated solutions of a system
\begin{equation}\label{system}
F(X)\ =\ 0
\end{equation}
where $F=(f_1,\ldots,f_n)$ are polynomials in the variables
$X=(x_1,\ldots,x_N)$.
First, a {\em homotopy} $H(X,t)$ is found with the following properties:
\begin{enumerate}
\item $H(X,1)= F(X)$.
\item The isolated solutions of $H(X,0)=0$ are known.
\item The system $H(X,t)=0$ defines finitely many (rational) curves 
$\sigma_i(t)$, and each isolated solution of~(\ref{system}) is connected to 
an isolated solution $\sigma_i(0)$ of  $H(X,0)=0$ by one of these curves.
\end{enumerate}
Given such a homotopy, numerical path continuation is used to trace these
curves from solutions of $H(X,0)=0$ to solutions of the original
system~(\ref{system}). 
When there are fewer solutions to $F(X)=0$ than to 
$H(X,0)=0$, some curves will diverge or become singular as
$t\rightarrow 1$, and it is expensive to trace such a curve.

When $N=n$, the system~(\ref{system}) is {\em square} and the homotopy
\begin{equation}\label{convex}
H(X,t)\quad :=\quad tF(X)\ +\ 
(1-t)G(X),
\end{equation}
where $G(X)=(x_1^{d_1}-a_1,\ldots,x_N^{d_N}-a_N)$
with $d_i:=\deg(f_i)$ and $a_i\neq 0$, gives
$\prod d_i$ curves.
This is the B\'ezout bound for a generic dense system $F$.

In practice, $F(X)=0$ may have fewer than 
$\prod d_i$ solutions and we desire a homotopy with no divergent curves.
Methods for such deficient systems which reduce the number of
divergent curves are
developed in~\cite{LSY_deficient,LS_deficient,LW_deficient}.
When the polynomials $f_1,\ldots,f_n$ have special
forms~\cite{MS_m-homogeneous,MSW_product}, then such 
homotopies~(\ref{convex}) are constructed where $G(X)$ shares
this special form.
When the polynomials $f_1,\ldots,f_n$ are sparse, polyhedral
methods~\cite{CVVerschelde,Huber_Sturmfels} give a homotopy.
The SAGBI homotopy algorithm (\S 2.3) is in the same spirit.
We exploit a special feature of the coordinate ring of the Grassmannian to
obtain a homotopy between the system (2.2) we wish to solve 
and one (2.12) whose solution may be obtained using polyhedral
methods.
Moreover, there are generically no divergent curves to be followed.

The overdetermined situation of $n>N$ is more delicate.
One difficulty is finding a homotopy $H(X,t)$ for an overdetermined system
as generic perturbations of $F$ have no solutions.
In~\cite[\S 2]{Sommese_Wampler_NAG},  this difficulty is avoided as
follows:
The  system $F=(f_1,\ldots,f_n)$ is replaced by $N$ random
linear combinations of  the $f_1,\ldots,f_n$ yielding a square system
whose isolated solutions include all isolated solutions of $F(X)=0$,
but typically many more.
They then find all isolated solutions of this random square subsystem.

For the Gr\"obner and Pieri homotopy algorithms, we gave (in \S\S 2.2
and 3.2--3) homotopies  $H(X,t)=(h_1(X,t),\ldots,h_n(X,t))$  and
solutions  $\sigma_i(0)$ at $t=0$ as above.  
For these, there are generically no divergent curves.
To efficiently follow the curves $\sigma_i(t)$, we select a square
subsystem  $MH(X,t)$ of $H(X,t)$ ($M$ is an $(N\times n)$-matrix).
If the Jacobian of $MH$ at each $\sigma_i(0)$ has the same rank ($N$) as
does the Jacobian of $H(X,t)$, then the curves $\sigma_i(t)$ remain 
components of the algebraic set defined by the equations
$$
MH(X,t)\ =\ 0.
$$
Moreover, other components of this set meet the curves $\sigma_i(t)$ 
in at most finitely many points $t$ in ${\bf C}-\{0\}$.
Thus,  we may use the square subsystem
$MH(X,t)$ to trace the curves  $\sigma_i(t)$ along some path 
from 0 to 1 in the complex plane.
We remark that in practice, $M$ may be chosen at random.

\section{Applications}
The algorithms of Sections 2 and 3 are useful for studying both
the pole assignment problem in systems theory~\cite{Byrnes} and real
enumerative geometry~\cite{Sottile_santa_cruz}.

We describe the connection to the control of linear systems
following~\cite{Byrnes}. 
Suppose we have a system (for example, a mechanical linkage) with inputs
$u\in {\bf R}^m$ and outputs $y\in {\bf R}^p$ for which there are internal
states $x\in {\bf R}^n$ such that the evolution of the system is governed by
the first order linear differential equation
\begin{equation}\label{linsystem}
\begin{array}{rcl}
\dot{x}&=&Ax + Bu,\\ 
    y&=&Cx.
\end{array}
\end{equation}
If the input is controlled by constant output feedback, $u=Fy$, then we obtain
$$
\dot{x}\ =\ (A+BFC)x.
$$
The natural frequencies of the controlled system are the roots 
$s_1,\ldots,s_n$ of 
\begin{equation}\label{charpoly}
\varphi(s)\ :=\  \det(sI-A-BFC).
\end{equation}
The pole assignment problem asks, given a system~(\ref{linsystem}) and a
polynomial $\varphi(s)$ of degree $n$, which feedback laws $F$ 
satisfy~(\ref{charpoly})?

A standard transformation~({\em cf.}~\cite[\S 2]{Byrnes}) transforms 
the input data $A,B,C$ into matrices
$N(s)$, $D(s)$ of polynomials with $\det(D(s))= \det(sI-A)$ and 
$N(s)D(s)^{-1}=C(sI-A)^{-1}B$
such that
\begin{equation}\label{schubert_form}
\varphi(s)\ =\ \det\left[
\begin{array}{cc} F & D(s)\\I& N(s)\end{array}\right].
\end{equation}
Here $I$ is the $p\times p$-identity matrix and the feedback law $F$ is an
$m\times p$-matrix.
If we let 
$$			
X\ :=\ 
\left[\begin{array}{c} F\\I\end{array}\right]
\qquad {\rm and}\qquad
K(s)\ :=\ 
\left[\begin{array}{c} D(s)\\N(s)\end{array}\right],
$$			
then $F$ gives local coordinates on $\mbox{\em Grass}(p,m+p)$
and~(\ref{schubert_form}) is equivalent to 
$$
X\cap K(s_i)\ \neq \ \{0\}\quad\mbox{for}\quad i\ =\ 1,\ldots,n.
$$
These conditions are independent for generic $A,B,C$  and
distinct $s_i$, hence
$n\leq mp$ is necessary for there to be any feedback laws $F$.
The critical case of $n=mp$ is an instance of the situation in \S 2.

In~\cite{Byrnes_Stevens_homotopy} homotopy
continuation was used to solve a specific feedback problem when 
$(m,p)=(3,2)$.
From this result, they deduced that the pole assignment
problem is not in general solvable by radicals. 
Despite this success, only few non-trivial examples have been computed
in the control theory literature~\cite{Rosenthal_Sottile}.

An important question is whether a given system may be controlled by
{\em real} output
feedback~\cite{Willems_Hesselink,Byrnes_real,Rosenthal_Schumacher_Willems}.
That is, if all roots of $\varphi(s)$ are real, are there real
feedback laws $F$ satisfying~(\ref{schubert_form})?  
Real enumerative geometry~\cite{Sottile_santa_cruz} asks a similar
question: 
Are there real linear subspaces $K_1,\ldots,K_n$ in general position 
with $\dim K_i=m+1-k_i$ and $k_1+\cdots+k_n=mp$ such that {\em all}
$p$-planes meeting each $K_i$ nontrivially are real?
When either $m$ or $p$ is 2~\cite{Sottile_real_lines}, 
$n\leq 5$~\cite{Sottile_explicit_pieri}, or when 
$m=p=3$ and the $k_i=1$~\cite{Sottile_santa_cruz}, the answer is yes.
In fact, the Pieri homotopies arose from these investigations.

B.~Shapiro and M.~Shapiro give a precise conjecture 
relating both applications.
Suppose 
$$
K_i(s)\ :=\ [\gamma(s),\gamma'(s),\gamma''(s),\ldots,\gamma^{(m+1-k_i)}(s)],
$$
where $\gamma(s)$ is a parameterization of a rational normal
(non-degenerate) curve in  ${\bf R}^{m+p}$ of degree $m+p-1$.
One such choice is
\begin{equation}\label{ratnorm}
\gamma(s)\ =\ {\it transpose}[1,s,s^2,\ldots,s^{m+p-1}].
\end{equation}
Geometrically, $K_i(s)$ is the $(m+1-k_i)$-plane which osculates the curve
$\gamma(s)$ at $s$.
Such osculating $m$-planes have been used to prove non-degeneracy results
in control theory. 

\begin{conj}[B.~Shapiro and M.~Shapiro]\label{conj:SS}
Let $s_1,\ldots,s_n$ be distinct real numbers and suppose 
$K_i(s_i)$ osculates $\gamma$ at $s_i$ and $k_1+\cdots+k_n=mp$.
Then each of the finitely many $p$-planes $X\subset {\bf C}^{m+p}$
satisfying  $X\cap K_i(s_i)\neq \{0\}$ for $i=1,\ldots,n$
is defined over the reals.
\end{conj}

\section{Computational results}
Our algorithms have been tested 
successfully in MATLAB, finding all 14 solutions in the
case $(m,p)=(4,2)$ for both the SABGI and Gr\"obner homotopy algorithm,
and all 15 solutions when $(m,p)=(6,2)$ and $k_1=\cdots=k_6=2$ for the
Pieri homotopy algorithm.

At present, the SAGBI and Gr\"obner homotopy algorithms have been fully
implemented.
Some timings from trial runs of these algorithms on a Sparc 20 
are displayed in Table~\ref{Sagbi_times}.
The input for these were $mp$ random complex $(m+p)\times m$-matrices.

We provide a comparison to methods based upon Gr\"obner bases.
Table~\ref{Sagbi_times} also gives the time on the  Sparc 20 for
the system Singular~\cite{SINGULAR} to compute a degree reverse
lexicographic Gr\"obner basis for the polynomial systems:
$$
\det \left[ X\mid K_i\right]\ =\ 0, \quad\mbox{for}\quad
i=3,\ldots,K_{mp}.
$$
Here $X$ is expressed in local coordinates for 
$\,\Omega_{[13]} \,\cap \,\Omega'_{[13]} \,$ and 
the $K_3,\ldots,K_{mp}$  are 
$(m+p)\times m$-matrices with random integral entries between $-4$ and
$4$. A degree reverse lexicographic Gr\"obner basis is the input for some
alternative numerical polynomial systems solvers (e.g.~eigenvalue
methods~\cite{Auzinger_Stetter}). 
We note that the Gr\"obner basis calculation 
did not terminate within one week in the case $(m,p)=(6,2)$.

\begin{table}[htb]
\begin{tabular}{cccccc}\hline\hline
$m$ & $p$&$d(m,p)$&SABGI homotopy&Gr\"obner homotopy&Gr\"obner basis\\ \hline
 3  &  2 &    5   &       $ {<}1$ & $    <0.5  $&$   <0.5$\\
 4  &  2 &   14   &       $   47$ & $     6    $&$    19  $\\
 5  &  2 &   42   &       $  373$ & $   408    $&$ 149,897$\\
 6  &  2 &  132   &       $3,364$ & $ 8,626    $&$\infty$\\ \hline\hline
\end{tabular}
\caption{Time (in seconds) \label{Sagbi_times}}
\end{table}

The final version of this paper will include data from implementations
of the Pieri homotopy algorithm, as well.

\end{document}